\documentclass[11pt]{amsart}
\usepackage{amssymb}
\newtheorem{theorem}{Theorem}[section]

\newtheorem{proposition}[theorem]{Proposition}

\newtheorem{lemma}[theorem]{Lemma}

\theoremstyle{remark}
\newtheorem{remark}[theorem]{Remark}

\newtheorem{example}[theorem]{Example}

\renewcommand{\O}{\on{O}}
\renewcommand{\o}{\on{o}}

\newcommand{\R}{\mathbb{R}}
\newcommand{\C}{\mathbb{C}}

\newcommand{\Z}{\mathbb{Z}}
\newcommand{\Q}{\ca{Q}}

\newcommand{\Cl}{{\on{Cl}}}
\newcommand{\Pin}{\on{Pin}}
\newcommand{\tpi}{2\pi\sqrt{-1}}
\newcommand\lie[1]{\mathfrak{#1}}
\renewcommand{\k}{\lie{k}}

\newcommand{\g}{\lie{g}}

\renewcommand{\t}{\lie{t}}

\newcommand{\on}{\operatorname}
 \newcommand{\Aut}{ \on{Aut} }

\newcommand{\ad}{ \on{ad} }

\newcommand{\End}{ \on{End} }

 \newcommand{\Spin}{ \on{Spin}}
 
\newcommand{\SO}{ \on{SO}} 
\renewcommand{\S}{\mathcal{S}}

\newcommand{\GL}{\on{GL}}
 \newcommand{\Vol}{ \on{Vol}}

\newcommand\dirac{/\kern-1.2ex\partial} 
\newcommand\qu{/\kern-.7ex/} 

\newcommand{\hra}{\hookrightarrow}

\renewcommand{\d}{{\mbox{d}}}
\newcommand{\ol}{\overline}

\newcommand\eps{\epsilon}
\newcommand\Om{\Omega}

\newcommand{\f}{\frac}

\newcommand{\p}{\partial}
\renewcommand{\l}{\langle}
\renewcommand{\r}{\rangle}
\newcommand{\hh}{{\f{1}{2}}}
\newcommand{\ti}{\tilde}

\newcommand\beqn{\begin{equation}}      
\newcommand\eeqn{\end{equation}}      
\newcommand{\ca}{\mathcal}
\newcommand{\wh}{\widehat}

\newcommand{\mf}{\mathfrak}
\newcommand{\beq}{\begin{eqnarray*}}
\newcommand{\eeq}{\end{eqnarray*}}

\newcommand{\id}{\on{id}}

\begin{document}

\title[ ]{Clifford algebras
and the \\ classical dynamical Yang-Baxter equation}

\author{A. Alekseev}
\address{University of Geneva, Section of Mathematics,
2-4 rue du Li\`evre, 1211 Gen\`eve 24, Switzerland}
\email{alekseev@math.unige.ch}

\author{E. Meinrenken}
\address{University of Toronto, Department of Mathematics,
100 St George Street, Toronto, Ontario M5S3G3, Canada }
\email{mein@math.toronto.edu}


\begin{abstract}
We describe a relationship of the classical dynamical Yang-Baxter equation 
with the following elementary problem for Clifford algebras: Given a 
vector space $V$ with quadratic form $\ca{Q}_V$, how is the exponential 
of an element in $\wedge^2(V)$ under exterior algebra multiplication 
related to its exponential under Clifford multiplication? 
\vskip.1in
\end{abstract}


\maketitle

\vskip 0.3cm

\section{Introduction}
Let $\g$ be a real Lie algebra, equipped with a non-degenerate
invariant quadratic form $\ca{Q}$. Let $\Theta\in\wedge^3\g$ be the
cubic element defined by the quadratic form and the Lie algebra
structure. An element $\mf{r}\in\wedge^2\g$ is called a {\em classical
$r$-matrix} for $\g$ if it satisfies the (modified) {\em classical Yang-Baxter
equation} (CYBE)
$$\hh [\mf{r},\mf{r}]^\g=\eps \Theta$$
for some {\em coupling constant} $\eps\in\R$. Here $[\mf{r},\mf{r}]^\g$ 
is defined using the extension of the Lie bracket 
to the Schouten bracket on the exterior algebra, 
$[\cdot,\cdot]^\g:\,\wedge^k\g\times\wedge^l\g\to
\wedge^{k+l-1}\g$.  Drinfeld \cite{dr:ha} and 
Semenov-Tian-Shansky \cite{se:wh} gave a geometric 
interpretation of the CYBE in terms of Poisson-Lie group 
structures, and a classification of $r$-matrices for semi-simple 
Lie algebras was obtained by Belavin-Drinfeld \cite{be:tr}. 

In 1994, G. Felder \cite{fe:con} described a generalization of the
CYBE called the classical {\em dynamical} Yang-Baxter equation
(CDYBE). Let $\k\subset \g$ be a Lie subalgebra. A classical dynamical
r-matrix is a $\k$-equivariant (meromorphic) function
$\mf{r}:\,\k^*\to\wedge^2\g$ satisfying the (modified) CDYBE
$$ \sum_i \f{\p \mf{r}}{\p \mu_i}\wedge e_i+\hh [\mf{r},\mf{r}]^\g=
\eps \Theta.$$
Here $e_i$ is a basis on $\k$ with dual basis $e^i\in\k^*$, and
$\mu_i$ are the corresponding coordinates on $\k^*$. Etingof-Varchenko
\cite{et:ge} interpreted classical dynamical r-matrices in terms of
Poisson Lie groupoids, and gave a classification for $\g$ semi-simple
and $\k$ of maximal rank. The classification was extended by
Schiffmann \cite{sc:on} to more general subalgebras, interpolating
between the Belavin-Drinfeld ($\k=0$) and the Etingof-Varchenko
cases. For more general $\g$, Etingof-Schiffmann \cite{et:mo}
introduced the moduli space of classical dynamical $r$-matrices, and
described its structure.

The CDYBE was discovered in conformal field theory and the theory of
quantum groups, but arises in many other contexts as well.
Balog-Feh\'{e}r-Palla \cite{ba:ch1,ba:ch2} discussed dynamical
$r$-matrices arising in WZNW theory and Feh\'{e}r-G\'{a}bor-Pusztai
\cite{fe:on} describe their appearance in Dirac reduction.  Poisson
geometric applications include work by Jiang-Hua Lu \cite{lu:cl}, who
proved that Poisson homogeneous structures on $G/T$ (for a compact Lie
group $G$ with maximal torus $T$) all come from solution of the CDYBE,
and Ping Xu \cite{xu:cl} who showed that solutions of the CDYBE for
$(\g,\k)$ give rise to certain invariant Poisson structures on
$\k^*\times G$.  The CDYBE also appears in the theory of quasi-Poisson
manifolds \cite{al:qu} and is related to the Duflo map in Lie theory
\cite{al:no}. Petracci \cite{pe:th} has studied the CDYBE for
arbitrary Lie algebras, and found relations to Lie's third theorem and
the Poincar\'{e}-Birkhoff-Witt isomorphism.

Our goal in this paper is to explain a close relationship of the CDYBE
with the following elementary problem in the theory of Clifford
algebras. Suppose $V$ is a vector space with a quadratic form $\Q_V$,
and let $\Cl(V)$ be its Clifford algebra. Recall that there is a
vector space isomorphism $q:\,\wedge V\to \Cl(V)$, called the
quantization map, and that elements $q(\wedge^2 V)\subset\Cl(V)$
exponentiate to the Spin group $\Spin(V)\subset \Cl(V)$. One can then
ask: How is the exponential $\exp(q(\lambda))\in\Cl(V)$ of
$\lambda\in\wedge^2(V)$ related to the corresponding exponential
$\exp(\lambda)\in\wedge(V)$ in the exterior algebra? Clifford calculus
gives nice formulas for this and closely related problems. Taking
$V=\g$, the answer to this problem produces solutions of the CDYBE.

\section{Clifford exponentials}\label{sec:cliffordexp}
In this section we describe various formulas for exponentials 
of quadratic elements in a Clifford algebra. Proofs of these 
formulas will be given in Section \ref{sec:proofs} below.

Let $V$ be a finite-dimensional real vector space, equipped with a
non-degenerate quadratic form $\ca{Q}_V$. The pair $(V,\ca{Q}_V)$ 
will be called a {\em quadratic vector space}. The Clifford algebra 
$\Cl(V)$ is the quotient of the tensor algebra $\ca{T}(V)$ by the 
ideal generated by elements $v\otimes v-\hh \Q(v),\ v\in V$. 
 
The involutive automorphism $a\in \Aut(\Cl(V))$ 
given on generators by $a(v)=-v$ gives $\Cl(V)$ the structure 
of a $\Z_2$-graded algebra. For the rest of this paper, 
commutators in the Clifford algebra, tensor products with 
other $\Z_2$-graded algebras and so on will always be taken in 
the $\Z_2$-graded sense. Let
$$q:\,\wedge V\to \Cl(V)$$ 
be the {\em quantization map}, defined by the inclusion $\wedge V\to
\ca{T}(V)$ as anti-symmetric tensors followed by the quotient map
$\ca{T}(V)\to \Cl(V)$. The quantization map $q$ is an isomorphism of
vector spaces, with inverse $q^{-1}$ the {\em symbol map}.

Let $\O(V)$ denote the orthogonal group of $(V,\Q_V)$, $\o(V)$ its Lie
algebra, and
$$  \lambda:\, \o(V)\to\wedge^2 V,\ \ 
\lambda(A)=\hh \sum_a A(e_a)\wedge e^a$$
the canonical isomorphism. Here $e_a$ is a basis of $V$ 
with dual basis $e^a\in V^*$ (identified with $V$ via $\Q_V$).
The map 
$$ \gamma:\,\o(V)\to \Cl(V),\ \ \gamma(A)=q(\lambda(A))$$
is a Lie algebra homomorphism into the even part of the Clifford algebra
(with bracket the commutator).  
$\gamma(A)$ generates the action of $A$ as a derivation of $\Cl(V)$, 
that is, 
$$A(v)=[\gamma(A),v]$$
for all $v\in V\subset \Cl(V)$. We will be interested in formulas 
for the Clifford algebra exponential $\exp(\gamma(A))$. One such 
formula reads (cf. \cite[Proposition 3.13]{be:he})
\begin{equation}\label{eq:c1} q^{-1}(\exp(\gamma(A)))={\det}^{1/2}(\cosh(A/2))
\exp(2 \lambda(\on{tanh}(A/2))),\end{equation}
where the square root of the determinant is a well-defined 
analytic function of $A$, equal to $1$ at $A=0$.  
If $\dim V$ is even, one has an alternative expression
\begin{equation}\label{eq:c2}
q^{-1}(\exp(\gamma(A)))=
{\det}^{1/2}(2\on{sinh}(A/2))
\exp(\iota_{\hh \lambda(\on{coth}(A/2)}))  \ \d\Vol
\end{equation}
Here $\d\Vol$ is the Euclidean volume form on $V$, for a given 
choice of orientation, and the square root of the determinant 
is defined as a Pfaffian. More generally, given commuting elements 
$c\in \O(V)$ and $A\in \o(V)$, and a lift  
$\hat{c}\in \Pin(V)$ one has the formula 
\begin{equation}\label{eq:c3}
q^{-1}(\hat{c}\exp(\gamma(A)))=\pm\,\, 
{\det}^{1/2}(c\exp(A)-I)
\exp\big(\iota_{\hh \lambda({\f{c\exp(A)+I}{c\exp(A)-I}})}\big) 
\ \d\Vol,
\end{equation}
where the sign depends on the choice of lift.

Most important for our purposes will be a formula relating
$\exp(\gamma(A))$ to the corresponding exponential $\exp(\lambda(A))$
in the exterior algebra.  Consider the holomorphic function
\begin{equation}\label{eq:j} 
j(z)=\f{\sinh(z/2)}{z/2}
\end{equation}
and let $f(z)$ be its logarithmic derivative, 
\begin{equation}\label{eq:f}
f(z)=(\ln j)'(z)=
 \frac{1}{2} \coth(\f{z}{2}) - \f{1}{z}.
\end{equation}
Note that $j$ is symmetric with simple zeroes at points 
$z\in 2\pi \sqrt{-1}\Z\backslash\{0\}$, while $f$ is anti-symmetric 
with simple poles at those points. The function $J\in C^\infty(\o(V))$
given by 
$$ J(A)=\det(j(A))$$
admits a unique smooth square root equal to $1$ at $A=0$. Define a 
meromorphic function $\mf{r}:\,\o(V)\to \wedge^2(V)$ by 
$$ \mf{r}(A)=\lambda(f(A))$$
and set $\S:\,\o(V)\to \wedge^{\on{even}}(V)$, 
$$ \S(A)=J^{1/2}(A)\exp(\mf{r}(A)).$$
\begin{theorem}\label{th:factor}
The function $\S$ is analytic on all of $\o(V)$.  Let $E$ be a vector
space of ``parameters'', and $\phi:\,V\to E$ a linear map with
components $\phi^a=\phi(e^a)$.  For all $A\in\o(V)$, the following
identity holds in $\Cl(V)\otimes\wedge(E)$:
\begin{equation}\label{eq:key}
 q\circ \iota_{\S(A)}\exp(\lambda(A)-{\sum}_a e_a \phi^a)
=
\exp(\gamma(A)-{\sum}_a e_a\phi^a).
\end{equation}
\end{theorem}

Notice that if $A$ has no eigenvalue equal to $0$ so that
$\lambda(A)\in\wedge^2(V)$ is non-degenerate, any element $\alpha\in
\wedge(V)$ can be written in the form $\alpha=\iota_\beta
e^{\lambda(A)}$ for a unique element $\beta\in \wedge(V)$. 
This explains formula \eqref{eq:key} for $\phi=0$ and 
$\lambda(A)$ non-degenerate. The remarkable feature of this formula 
is that $\S$ extends analytically to all $A$, and is independent of  
$\phi$. 

For applications to Lie algebras, it is sometimes useful to 
write the right hand side of \eqref{eq:key} somewhat differently. 
Define holomorphic functions 
$$ 
g(z)=\f{\sinh(z)-z}{z^2},\ \ 
j^R(z)=\f{e^z-1}{z},\ \ j^L(z)=\f{1-e^{-z}}{z}.
$$
Let $\psi=\phi\circ j^R(A):\,V\to E$, with components
$\psi^a=\psi(e^a)$, and let $\varpi(A)\in\wedge^2(E)$ be the image of
$\lambda(g(A))\in\wedge^2(V)$ under the extended map
$\phi:\,\wedge(V)\to \wedge(E)$.  Then
\begin{equation}\label{eq:alter}
q\circ \iota_{\S(A)}\exp(\lambda(A)-{\sum}_a e_a \phi^a)
=\exp(-\varpi(A)) \exp(\gamma(A))\exp(-\sum_a e_a \psi^a).
\end{equation}

\section{Quadratic Lie algebras}\label{sec:quadr}
In this Section, we specialize Theorem \ref{th:factor}
to Lie algebras with an invariant quadratic form. We show that  
the exponentials in the Clifford and exterior algebras 
satisfy natural differential equations. In the following 
Section the dynamical Yang-Baxter equation emerges as a consistency 
condition for these differential equations.

\subsection{Lie algebra consequences of Theorem \ref{th:factor}}
A {\em quadratic Lie algebra} is a Lie algebra $\g$, together with an
invariant, non-degenerate quadratic form $\Q$. We will denote 
the Lie bracket by $[\cdot,\cdot]^\g$, to avoid confusion with 
commutators. The invariance condition 
means that the adjoint representation $\ad:\,\g\to \End(\g)$ takes
values in $\o(\g)$. Examples of quadratic Lie algebras include semi-simple
Lie algebras, and semi-direct products $\g=\mf{s}\ltimes
\mf{s}^*$, where $\mf{s}$ is {\em any} real Lie algebra, acting on its
dual $\mf{s}^*$ by the coadjoint action. Also, given a possibly 
degenerate invariant quadratic form $\Q'$ on a Lie algebra $\g'$, 
the quotient of $\g=\g'/K$ by the radical $K$ of the quadratic form is a 
quadratic Lie algebra. Note that quadratic Lie algebras
are unimodular. See the work of Medina-Revoy \cite{me:al1} 
for further information and classification results.

The ingredients in Theorem \ref{th:factor}, and its consequences, take
on geometric meanings if $V=\g$ is a quadratic Lie algebra, and 
$A=\ad_\mu$ for $\mu\in\g$.  

\begin{enumerate}
\item[$\lambda(A),\gamma(A)$] 
The map $\lambda^\g=\lambda\circ \ad:\,\g\to \wedge^2\g$ is the 
map dual to the Lie bracket. The map $\gamma^\g=q(\lambda^\g):\,\g\to \Cl(\g)$ 
generates the adjoint action of $\g$ on the Clifford algebra.
\vskip.03in
\item[$j^{L/R}(A)$]
Let $G$
be the connected, simply connected Lie group having $\g$ as its Lie
algebra, and let $\exp_G:\,\g\to G$ be the exponential map.  Denote by
$\theta^L,\theta^R\in\Om^1(G,\g)$ the left/right invariant
Maurer-Cartan forms. 
Under left trivialization of the tangent bundle $TG$, the differential
$\exp_G$ at $\mu\in \g$ is given by the operator $j^L(\ad_\mu)$
(cf. \cite[Theorem II.1.7]{he:di}). Equivalently, the 
value of $\exp^*\theta^L$ at $\mu$ is 
given by 
$$(\exp^*\theta^L)_\mu=j^L(\ad_\mu).$$
Similarly one has $(\exp^*\theta^R)_\mu=j^R(\ad_\mu)$. 
\item[$J(A)$] Let $J^\g=J\circ \ad$. The quadratic form $\Q$ defines
a translation invariant measure on $\g$ and a bi-invariant measure on
$G$. Since $\g$ is unimodular, 
$\det(j^L(\ad_\mu))=\det(j^R(\ad_\mu))=J^\g(\mu)$.  Thus $J^\g$
is the Jacobian of the exponential map $\exp_G$ with respect to left or 
right trivialization of $TG$, and the subset 
of $\g$ where $A=\ad_\mu$ has eigenvalues in $2\pi i \Z\backslash\{0\}$ 
is the set of critical points. 
\item[$\varpi(A)$] 
Let
$ B(\mu,\mu')= \hh(\Q(\mu+\mu')-\Q^\g(\mu)-\Q(\mu'))$
denote the 
symmetric bilinear form associated with $\Q$. Let
$\eta=\f{1}{12}B(\theta^L,[\theta^L,\theta^L]^\g)\in\Om^3(G)$ be the
{\em Cartan 3-form} on $G$. It is bi-invariant and therefore 
closed. Let $\varpi^\g \in \Om^2(\g)$ be the image of $\exp_G^*\eta$
under the usual homotopy operator $\Om^p(\g)\to \Om^{p-1}(\g)$, so
that $\d\varpi^\g=\exp_G^*\eta$. If we identify $\wedge^2T^*_\mu\g \cong \wedge^2\g$, the value of $\varpi^\g$ at $\mu$ is given by the formula
(cf. \cite{al:mom}),
$$ \varpi^\g_\mu=\varpi(\ad_\mu).$$ 
\item[$\mf{r}(A)$] Let $\mf{r}^\g(\mu)=\mf{r}(\ad_\mu)=\lambda(f(\ad_\mu))$. 
Given $\xi\in\g$ let $\xi^L,\xi^R$ denote the left/right invariant 
vector fields on $G$ generated by $\xi$. On the subset of $\g$ 
where $\exp_G$ is regular, the vector field 
$\hh\exp_G^*(\xi^L+\xi^R)$ is well-defined. It differs from 
the constant vector field $\xi$ by a vector field tangent to 
orbit directions for the adjoint action. It turns out \cite[Lemma A.1]{al:no}
that the difference at $\mu\in\g$ coincides with the vector 
field generated by $f(\ad_\mu)\xi$.  
\end{enumerate}

Theorem \ref{th:factor} takes on the following form.  
Let $\S^\g(\mu)=\S(\ad_\mu)$. 

\begin{proposition}\label{prop:identity}
The identity 
\begin{equation}\label{eq:key1}
q\circ \iota_{\S^\g} (e^{\lambda^\g -\sum_a e_a \phi^a})= 
e^{\gamma^\g-\sum_a e_a \phi^a}
\end{equation}
holds in $\Cl(\g)\otimes C^\infty(\g)\otimes\wedge E$.
\end{proposition}
The alternative formula \eqref{eq:alter} takes on 
a particularly nice form for $E=T^*_\mu\g$, with 
$\phi:\,\g\cong\g^*\to T^*_\mu\g$ the standard 
identification. Then $\phi^a=\d\mu^a$ 
where $\mu^a$ are the coordinate functions on $\g$, and 
$\psi^a=\exp_G^*(\theta^L)^a$. 
Since $G$ is simply connected, the Lie algebra homomorphism 
$\gamma^\g:\,\g\to \Cl(\g)$ exponentiates to a Lie group homomorphism 
$\tau:\,G\to \Spin(\g)$, with $e^{\gamma^\g}=\exp_G^*\tau$. The 
resulting formula 
$$q\circ \iota_{\S^\g} (e^{\lambda^\g 
-\sum_a e_a d \mu^a})=  e^{-\varpi^\g}  \exp_G^* 
\left(\tau\  e^{-\sum_a e_a (\theta^L)^a}\right)
$$
relates $\tau\  e^{-\sum_a e_a (\theta^L)^a}\in \Cl(\g)\otimes\Om(G)$ 
and $e^{\lambda^\g -\sum_a e_a d \mu^a}\in\wedge\g\otimes\Om(\g)$. 
For the case of compact Lie algebras, this result was proved in 
\cite[Section 6.3]{al:no}. 

\subsection{Lie algebra differential}
In this Section,  $\g$ denotes an arbitrary Lie algebra (not necessarily 
quadratic), and
$\lambda^\g:\,\g^*\to\wedge^2\g^*$ the map dual to the Lie bracket,
\begin{equation}\label{eq:lambda}
\iota_\xi\iota_\eta \lambda^\g(\mu)=\l\mu,[\xi,\eta]^\g\r.\end{equation}
In a basis $e_a$ of $\g$, with dual basis $e^a$ of $\g^*$, 
$${\textstyle  \lambda^\g(e^c)=-\hh \sum^c_{ab} f^c_{ab} e^a\wedge e^b}$$
where $f^c_{ab}=\l e^c,[e_a,e_b]^\g\r$ are the structure constants. 
Recall that the {\em Lie algebra differential} 
$\d^\g:\,\wedge^\bullet\g^*\to \wedge^{\bullet+1}\g^*$ 
is the (degree $+1$) derivation given 
on generators $\mu\in\g^*$ by
\begin{equation}\label{eq:dwedge}\d^\g\mu=\lambda^\g(\mu).\end{equation}
\begin{lemma}\label{lem:lambda}
Let $E$ be some vector space, $\phi:\,\g\to E$ a linear map, 
and $\phi_a=\phi(e_a)$.  Then 
$$\lambda^\g-\sum_a e^a\phi_a$$ 
is closed under the differential, $\d^\g+\sum_a \f{\p}{\p \mu_a}\phi_a$
on $C^\infty(\g^*)\otimes \wedge\g^*\otimes\wedge E$
. 
\end{lemma}
\begin{proof}
Since $\d^\g\lambda^\g(\mu)=\d^\g\d^\g\mu=0$, 
this follows from the calculation, 
$$\d^\g(\sum_a e^a \phi_a)= \sum_a \lambda^\g(e^a) \phi_a=\sum_a
\f{\p \lambda^\g}{\p \mu_a} \phi_a .$$
\end{proof}

We will need the following Lemma, 
describing the transformation 
of $\d^\g$ under conjugation with $\exp(\iota_{\mf{r}})$:
\begin{lemma}\label{lem:rconj}
For any $\mf{r}\in \wedge^2\g$, 
\begin{equation}\label{eq:conj}
\exp(-\iota_{\mf{r}}) \circ \d^\g\circ \exp(\iota_{\mf{r}})
=\d^\g-\hh \iota_{[\mf{r},\mf{r}]^\g}
+\sum_a e^a \circ \iota_{[e_a,\mf{r}]^\g}
-\iota_{\mf{u}}
\end{equation}
where $\mf{u}\in\g$ is the image of $\mf{r}$ under the Lie bracket map 
$\wedge^2\g\to \g,\ \xi\wedge \xi'\mapsto [\xi,\xi']^\g$. 
\end{lemma}
\begin{proof}
Write $\iota_a=\iota_{e_a}$.
Then $\d^\g=-\hh \sum_{abc} f_{ab}^c e^a \wedge e^b \circ \iota_c$.
Introduce 
components $\mf{r}^{ab}$ by $ \mf{r}=\hh\sum_{ab} \mf{r}^{ab}
e_a\wedge e_b$. Then the Schouten bracket of $\mf{r}$ with itself is 
given by the formula, 
$$
[\mf{r},\mf{r}]^\g= \sum_{abc}(\sum_{kl} \mf{r}^{ak}f_{kl}^b
 \mf{r}^{lc}) e_a\wedge e_b \wedge e_c .$$
We compute the left hand side of \eqref{eq:conj} as
a sum $\sum_{j=0}^\infty \f{1}{j!}\ad^j(-\iota_{\mf{r}})\d^\g$: 
\beq \ad(-\iota_{\mf{r}})\d^\g&=&
\f{1}{4}\sum_{abklm}
\mf{r}^{ab}f^k_{lm}[\iota_a\iota_b,e^l \wedge e^m]\iota_k=-
\sum_{abkm}
\mf{r}^{ab}f_{bm}^k
e^m\iota_a\iota_k-\hh \sum_{abk} \mf{r}^{ab} f_{ab}^k \iota_k \\
&=&\sum_a e^a \circ \iota_{[e_a,\mf{r}]^\g} -\iota_{\mf{u}}
,\\
 \ad^2(-\iota_{\mf{r}})\d^\g&=&
\hh \sum_{abkmst}
\mf{r}^{st}\mf{r}^{ab}f_{bm}^k [\iota_s\iota_t, e^m]\iota_a\iota_k=
\sum_{abkst}
\mf{r}^{ab}\mf{r}^{st}f_{bt}^k  \iota_s\iota_a\iota_k\\
&=&-\iota_{[\mf{r},\mf{r}]^\g}\\
\ad^m(-\iota_{\mf{r}})\d^\g&=&0,\ \ \ \ m\ge 3
\eeq
\end{proof}

\subsection{Clifford differential}
Let $(\g,\ca{Q})$ be a quadratic Lie algebra.
Similar to \eqref{eq:dwedge}
there is a unique odd derivation 
$\delta^\g$ on $\Cl(\g)$ given on generators $\mu\in\g$ by 
$$ \delta^\g\mu=\gamma^\g(\mu).$$
In fact $\delta^\g$ may be
written as a ($\Z_2$-graded) commutator: Let $\Theta\in\wedge^3\g$ be
the cubic element defined by $\iota_\mu \Theta=\lambda^\g(\mu).$ In
terms of a basis $e_a$ of $\g$, with dual basis $e^a$, we have
$$ \textstyle{ \Theta=-\f{1}{6}\sum_{abc} f^{abc} e_a \wedge e_b \wedge e_c}$$
where $f^{abc}=B(e^a,[e^b,e^c]^\g)$.  Then
$\delta^\g=[q(\Theta),\cdot]$. Kostant-Sternberg \cite{ko:sy2}
made the beautiful observation that $q(\Theta)$ squares to a
constant, hence that $\delta^\g$ squares to $0$. We will call
$\delta^\g$ the {\em Clifford differential}.  Under the
quantization map $q$, the Lie algebra and Clifford differentials are
related as follows \cite[Proposition 3.3]{al:no}:
\begin{equation}\label{eq:related}
{\textstyle q^{-1}\circ \delta^\g\circ q =\d^\g+\f{1}{4}\iota_{\Theta}}.
\end{equation}
Replacing $\lambda$ with $\gamma$ in the proof of Lemma \ref{lem:lambda}, 
we find: 
\begin{lemma}\label{lem:gamma}
Let $E$ be some vector space, $\phi:\,\g\to E$ a linear map, 
and $\phi^a=\phi(e^a)$.  Then 
$$\gamma^\g-\sum_a e_a\phi^a$$
is closed under the differential, 
$\delta^\g+\sum_a \f{\p}{\p \mu^a}\phi^a$ on 
$C^\infty(\g)\otimes \Cl(\g)\otimes\wedge E$. 
\end{lemma}

\section{Solutions of the classical dynamical Yang-Baxter equation}
In this section, we will use our Clifford algebra techniques to 
construct solutions to the CDYBE in a number of cases. 
\begin{theorem}\label{th:cdybe} 
The function $\mf{r}^\g(\mu)=\lambda(f(\ad_\mu))$ for $f(z)=
\hh\coth(z/2)-\f{1}{z}$ solves the CDYBE for $\k=\g$, 
with coupling constant $\eps=\f{1}{4}$.
\end{theorem}
\begin{proof}
The proof relies on the identity \eqref{eq:key1} from 
Proposition \ref{prop:identity}. 
To simplify notation, we denote 
$\mf{r}^\g,\ca{S}^\g$ simply by $\mf{r},\ca{S}$, respectively.

By Lemma \ref{lem:gamma}, the right hand side of \eqref{eq:key1} is closed 
under the differential, $\delta^\g+\sum_a \f{\p}{\p \mu^a}\phi^a$. 
Hence 
$$ (\delta^\g+\sum_a \f{\p}{\p \mu^a}\phi^a)
\left(q\circ \iota_{\S} (e^{\lambda^\g -\sum_a e_a \phi^a})\right)
=0.$$
By \eqref{eq:related}, this gives 
$$ (\d^\g+\f{1}{4}\iota_{\Theta}+\sum_a \f{\p}{\p \mu^a}\phi^a)
\left(\iota_{\S} (e^{\lambda^\g -\sum_a e_a \phi^a})\right)
=0.$$

Since $\iota_\S = J^{1/2}\exp(\iota_{\mf{r}})$ we have, using Lemma
\ref{lem:rconj}, 
\beq 
\d^\g\circ \iota_\S&=&\iota_\S\circ\Big(
\d^\g-\iota_{\mf{u}}+
\sum_a e^a\circ \iota_{[e_a,\mf{r}]^\g}
-\hh \iota_{[\mf{r},\mf{r}]^\g}
\Big)\eeq
Furthermore,
\beq
  \f{\p }{\p \mu^a}  \circ \iota_\S&=& 
\iota_\S \circ 
\Big( \f{\p }{\p \mu^a} +\hh   \f{\p \ln(J)}{\p \mu^a}
+\iota_{ \f{\p \mf{r}}{\p \mu^a}}\Big).
\eeq
Using that $e^{\lambda^\g- \sum_a e_a \phi^a}$ is closed under the
differential $\d^\g+\sum_a\f{\p }{\p \mu^a}\phi^a$, we 
therefore obtain
\begin{equation}\label{eq:star}
\Big(-\iota_{\mf{u}}+
\sum_a e^a\circ \iota_{[e_a,\mf{r}]^\g}
-\hh \iota_{[\mf{r},\mf{r}]^\g}
+\sum_a \left(\hh   
\f{\p \ln(J)}{\p \mu^a}
+\iota_{ \f{\p \mf{r}}{\p \mu^a}}\right)\phi^a
+\f{1}{4}\iota_{\Theta} \Big)
e^{\lambda^\g- \sum_b e_b \phi^b}=0.
\end{equation}
Multiply this Equation from the left by $\exp(\sum_a e_a
\phi^a)$, and pick the coefficient  {\em cubic} in $\phi$'s. 
Only
the three terms involving $\sum_a \phi^a\,\iota_{\f{\p \mf{r}}{\p \mu^a}}$, 
$\iota_{[\mf{r},\mf{r}]^\g}$ and  $\iota_{\Theta}$ contribute to this
coefficient, and we obtain: 
$$\sum_a  \phi\left(\f{\p \mf{r}}{\p \mu^a}\right)\wedge\phi^a 
+\hh \phi([\mf{r},\mf{r}]^\g)
-\f{1}{4}\phi(\Theta)\equiv
\phi\left(\f{\p \mf{r}}{\p \mu^a}\wedge e^a 
+\hh [\mf{r},\mf{r}]^\g
-\f{1}{4} \Theta\right)=0.
$$
Taking $E=\g$, with $\phi$ the identity map, this is exactly the 
CDYBE. 
\end{proof}

\begin{remark}
It is not hard to work out the 
coefficients of $\phi$ of lower degree. Two of these 
identities simply state that $\mf{r}$ is equivariant
and $J$ is invariant. The remaining identity reads
(cf. \cite[Lemma A.2]{al:no})
$$\hh \sum_a \f{\p \ln(J)}{\p \mu^a}e^a +\mf{u}=0.$$
\end{remark}

Theorem \ref{th:cdybe} was first obtained by Etingof-Varchenko \cite{et:ge}
in the semi-simple case. See Etingof-Schiffmann \cite{et:mo} for a 
proof in the quadratic case.

More generally, let 
$\k\subset \g$ be a {\em quadratic subalgebra} of $\g$, i.e. a subalgebra 
such that the restriction of $\Q$ to $\k$ is non-degenerate. 
Let $\mf{p}$ be the orthogonal complement of $\k$, so that 
$\g=\k\oplus \mf{p}$. Suppose that for $\mu$ in an open dense subset 
of $\k$, the operator $\ad_\mu^{\mf{p}}:=\ad_\mu|_{\mf{p}}$ is 
invertible. Let $\ad_\mu^\k=\ad_\mu|_\k$ and define $\mf{r}^{\k}:\,\k\to\wedge^2\k$ and 
$\mf{r}^{\mf{p}}:\,\k\to\wedge^2\mf{p}$ by 
$$ \mf{r}^{\k}(\mu)=\lambda(f(\ad_\mu^\k)),\ \ 
\mf{r}^{\mf{p}}(\mu)=\hh \lambda(\coth(\ad_\mu^{\mf{p}}/2)),\ \ 
\mf{r}=\mf{r}^\k+\mf{r}^{\mf{p}}
$$
\begin{theorem}
The function $\mf{r}=\mf{r}^\k+\mf{r}^{\mf{p}}$ solves the 
CDYBE for $\k\subset \g$, with coupling constant $\eps=\f{1}{4}$.
\end{theorem}

\begin{proof}
Let 
$$ J^{\k}(\mu)={\det}(j(\ad_\mu^\k)),\ \ 
J^{\mf{p}}(\mu)=\det(2\sinh(\ad_\mu^{\mf{p}}/2)),\ \ 
J=J^{\k}J^{\mf{p}},$$
and 
$$ \S^{\k}=(J^{\k})^{1/2}\,\exp(\mf{r}^{\k}),
\S^{\mf{p}}=(J^{\mf{p}})^{1/2}\,\exp(\mf{r}^{\mf{p}}),\ \ 
\S=\S^{\k}\S^{\mf{p}}.$$
Here the square root $(J^{\mf{p}})^{1/2}$ is defined as a Pfaffian,
for some choice of orientation on $\mf{p}$. Let $\d\Vol_{\mf{p}}$ be 
the volume form defined by the orientation and the quadratic form
$\Q|_{\mf{p}}$.

Let $e_a$ be a basis of $\g$ given by a basis of $\k$ followed by a
basis of $\mf{p}$. In what follows, summation over $a$ denotes
summation over the entire basis, while summation over $i$ denotes
summation over the basis of $\k$.

The restriction of $\gamma^\g$ to $\k$ is a sum
$\gamma^\g|_\k=\gamma^\k+\gamma^{\mf{p}}$, where $\gamma^{\mf{p}}$
takes values in $\Cl(\mf{p})$.  Combining Equation \eqref{eq:key1},
with $\k$ in place of $\g$, with Equation \eqref{eq:c2}, for
$V=\mf{p}$, we obtain the following identity in $C^\infty(\k)\otimes
\Cl(\g)\otimes \wedge E$:
\begin{equation}\label{eq:key2}
q\circ \iota_{\S} (e^{\lambda^\k -\sum_i e_i \phi^i}\wedge 
\d\Vol_\mf{p})= 
e^{\gamma^\g|_\k-\sum_i e_i \phi^i}.
\end{equation}
Write  $\d^\g=\d^\k+\d'$, where $\d^\k$ is 
extended to $\wedge\g$ by letting $\d^\k\mu=0$ for $\mu\in\mf{p}$.
Since $\gamma^\g|_\k-\sum_i e_i \phi^i$ is closed 
under the differential 
$\delta^\g+\sum_i \f{\p}{\p \mu^i} \phi^i$, and 
$\lambda^\k- \sum_i e_i \phi^i$ is closed under the 
differential $\d^\k+\sum_i \f{\p}{\p \mu^i} \phi^i$
we can proceed as 
in the proof of Theorem \ref{th:cdybe} to obtain 
$$ 
\Big(\d'-\iota_{\mf{u}}+
\sum_a e^a\circ \iota_{[e_a,\mf{r}]^\g}
-\hh \iota_{[\mf{r},\mf{r}]^\g}
+\sum_i\big(\hh   \f{\p \ln(J)}{\p \mu^i}
+\iota_{ \f{\p \mf{r}}{\p \mu^i}}\big)\phi^i
+\f{1}{4}\iota_{\Theta} \Big)
\big(e^{\lambda^\k- \sum_i e_i \phi^i}\wedge \d\Vol_{\mf{p}}\big)=0.
$$
Multiply from the left by $e^{\sum_a e_a\phi^a-\lambda^{\mf k}}$
(summation over the entire basis of $\g$). The term cubic in $\phi$'s
is proportional to $\d\Vol_{\mf{p}}$, and the coefficient gives the
CDYBE.  
\end{proof}

Still more generally, suppose $c\in \O(\g)$ is an automorphism of
$\g$.  Suppose $\k$ is a quadratic subalgebra contained in the fixed
point set of $c$. Let $\mf{p}$ be its orthogonal complement as above,
and suppose that for $\mu$ in an open dense subset of $\k$, the
operator $c\exp(\ad_\mu)-I$ is invertible on $\mf{p}$. Then
$$ J^{\mf{p}}_c={\det}_{\mf{p}}(c\exp(\ad_\mu^{\mf{p}})-I),\ \
\mf{r}^{\mf{p}}_c=\hh\lambda\Big(\f{c\exp(\ad_\mu)+I}{c\exp(\ad_\mu)-I}
\Big|_{\mf{p}}\Big),\ \ 
\S^{\mf{p}}_c=(J^{\mf{p}}_c)^{1/2}\,\exp(\mf{r}^{\mf{p}}_c)
$$
are well-defined meromorphic functions on $\k$.

\begin{theorem}
The function $\mf{r}=\mf{r}^\k+\mf{r}^{\mf{p}}_c$ solves the CDYBE 
with coupling constant $\eps=\f{1}{4}$.
\end{theorem}
\begin{proof}
Equations \eqref{eq:c3} and \eqref{eq:key1} give the 
following  identity in 
$C^\infty(\k)\otimes\Cl(\g)\otimes  \wedge E$:
\begin{equation}\label{eq:key3}
q\circ \iota_{\S} (e^{\lambda^\k -\sum_i e_i \phi^i}\wedge 
\d\Vol_\mf{p})= \pm\,\,
\hat{c}\, e^{\gamma^\g|_\k-\sum_i e_i \phi^i}
\end{equation}
where $\S=\S^\k\,\S^{\mf{p}}_c$. 
The element $\wh{c}$ commutes with $q(\Theta)$, since $c$ is an automorphism 
of $\g$ preserving the quadratic form. Hence $\wh{c}$ 
commutes with $\delta^\g$, and hence the right hand side of 
\eqref{eq:key3} is closed under $\delta^\g+\sum_i \f{\p}{\p \mu^i}\phi^i$. 
The rest of the proof is as before. 
\end{proof}

The classical dynamical r-matrix described here was first obtained by
Etingof-Schiffmann \cite{et:mo}, for the case that $c$ is a finite
order automorphism and $\k=\g_0$. (Note that the fixed point set of a
finite order automorphism $c\in O(\g)$
is a quadratic subalgebra.)

\begin{example}
Let $(\k,\Q^\k)$ be a quadratic Lie algebra, and $\g=\k^\C$ its 
complexification. The real part of the complexification of $\Q^\k$ 
defines a non-degenerate quadratic form $\Q$ on $\g$, with 
$\mf{p}=\sqrt{-1}\k$ (viewed as a real subalgebra). Let $c\in \O(\g)$ 
denote the automorphism given by complex conjugation. The $\mf{r}$-matrix 
described above has the form
$$ \mf{r}=\lambda\left(f(\ad_\mu^\k)+\hh \tanh(\ad_\mu^{\mf{p}}/2)\right).$$
\end{example}
\vskip.2in

As explained in \cite{et:ge,et:le}, one obtains other solutions of the CDYBE 
by scaling or taking limits: 

(i) If $\mf{r}$ is a solution of the CDYBE for coupling constant 
$\eps$, then $\mf{r}_t(\mu):=t^{-1}\ \mf{r}(t^{-1}\mu)$ is a solution with coupling 
constant $t^{-2}\eps$. Applying this to the trigonometric solutions 
obtained above, and taking the limit for $t\to \infty$ one obtains 
{\em rational} solutions of the CDYBE with vanishing coupling constant. 
If $\mf{r}$ is anti-symmetric in $\mu$, one can also take 
imaginary $t$ changing the sign of the coupling constant. This 
replaces $\coth$ with $\cot$ in our formulas. 

(ii) For any element $\nu$ in the center of $\k$, the shifted
$r$-matrix $\mu\mapsto \mf{r}(\mu+\nu)$ again solves the CDYBE for
coupling constant $\eps$. Furthermore, if the limit
$$ \mf{r}_\nu(\mu)=\lim_{t\to \infty}\mf{r}(\mu+t\nu)$$ 
exists, then the limiting $r$-matrix again solves the CDYBE for $\eps$. 
For instance, if $\g$ is semi-simple and $\k=\t$, one obtains 
constant r-matrices by taking $\nu\in\t$ some regular element. 

(iii) Recall that classical dynamical $r$-matrices 
$\mf{r}:\,\k\to \wedge^2\g$ are always required to be 
$\k$-invariant. Hence, the Schouten bracket of any 
element of $\wedge\k$ with $\mf{r}$ vanishes. Thus if 
$\mf{s}:\,\k\to \wedge^2\k$ solves the CDYBE for $(\k,\k)$ with 
coupling constant 
$\delta$, then $\mf{r}+\mf{s}$ solves the CDYBE for $(\g,\k)$ with 
coupling constant $\delta+\eps$. In particular, if $\k$ is 
Abelian, any closed 2-form on $\k$ gives rise to a solution 
of the CDYBE for $(\k,\k)$ with coupling constant $0$. 

\section{Clifford algebra calculations}\label{sec:proofs}
In this Section we prove the Clifford algebra identities 
from Section \ref{sec:cliffordexp}. 
These formulas are most systematically obtained from the {\em spinor
representation} for the Clifford algebra of the direct sum $V\oplus
V^*$, which we briefly review.

\subsection{Spinor representation}\label{subsec:spinor}
Let $V$ be a finite-dimensional real vector space. The 
direct sum  $W=V\oplus V^*$ carries a quadratic form 
$$\textstyle{ \Q_W(v\oplus\alpha)=2\,\alpha(v)}.$$
Let $\Cl(W)$ be the Clifford algebra of $(W,\Q_W)$, and 
consider the algebra representation 
$$ \pi:\,\Cl(W)\to \mf{gl}(\wedge V)$$
where generators $v\in V$ act by wedge product and generators 
$\alpha\in V^*$ act by contraction. The restriction of 
$\pi$ to a group representation of $\Spin(W)\subset 
\Cl(W)^\times$ is called the spinor representation. 

The group $\SO(W)$ contains $\wedge^2(V),\wedge^2(V^*),\GL(V)$ 
as distinguished subgroups, lifting to subgroups 
$\wedge^2(V),\wedge^2(V^*),\on{ML}(V)$ of $\Spin(W)$:

(i) For any skew-adjoint linear map $D:\,V^*\to V$ let 
$\lambda(D)=\hh \sum_a D(e^a)\wedge e_a\in\wedge^2(V)$. 
There is an inclusion 
$$\wedge^2(V)\to \SO(W),\ \ \lambda(D)\mapsto \left(\begin{array}{cc}
I&D\\0&I \end{array}\right).$$ 
This inclusion lifts to a vector 
subgroup $\wedge^2(V)\hra \Spin(W)$, and the action of 
$\pi(\lambda(D))$ is wedge product with 
$\exp(\lambda(D))$. 

(ii) Similarly, there is an identification $E\mapsto \lambda(E)$ 
of skew-adjoint linear maps $V\to V^*$ with $\wedge^2 V^*$. The  
inclusion 
$$\wedge^2 V^*\to \SO(W),\ \lambda(E)\mapsto 
\left(\begin{array}{cc}
I&0\\ E&I
\end{array}\right).$$
lifts to an inclusion $\wedge^2 V^*\hra \Spin(W)$, and  
$\pi(\lambda(E))$ is given by contraction with $\exp(\lambda(E))$. 

(iii) Finally, there is an inclusion 
$$ \GL(V)\to \SO(W),\ \ R\mapsto \left(\!\begin{array}{cc}R&0\\0&(R^{-1})^*
\end{array}\!\!\!\right).$$
The metalinear group $\on{ML}(V)\hra \Spin(W)$ is the inverse image of 
$\GL(V)$ under this map. 
The action of an element $\hat{R}\in \on{ML}(V)$, covering 
$R\in \GL(V)$,  in the spinor representation is given by 
\begin{equation}\label{eq:R}
\pi(\hat{R}).\alpha=\f{R.\alpha}{{|\det|}^{1/2}(\hat{R})}.
\end{equation}
Here ${|\det|}^{1/2}:\,\on{ML}(V)\to \R^\times$ is 
a suitable 
choice of square root of $|\det|:\,\on{GL}(V)\to \R_{>0}$ 
(defined by this 
formula), and $R.\alpha$ is defined 
by the unique extension of $R\in\GL(V)$ to an algebra automorphism 
of $\wedge(V)$.

\subsection{The action of $\Spin(V)$ on $\wedge V$}
We now return to our original setting, where $V$ itself comes 
equipped with a quadratic form $\ca{Q}_V$. 
Under the identification $q:\,\wedge V\cong \Cl(V)$, 
the left multiplication of the Clifford algebra on itself 
defines a representation, 
\begin{equation}\label{eq:pi}
\varrho:\,\Cl(V)\to \mf{gl}(\wedge V),
\end{equation}
given on generators $v\in V$ by  
$$\varrho(v).\alpha=v\wedge\alpha+\hh
\iota_v\alpha.$$
The symbol map $q^{-1}:\,\Cl(V)\to \wedge(V)$ can be 
expressed in terms of $\varrho$ as $q^{-1}(x)=\varrho(x).1,$
the action on $1\in\wedge V$. 

We will now relate $\varrho$ (hence also the symbol map) to the 
representation $\pi$ from \ref{subsec:spinor}. 
Let $\ol{V}$ denote the 
same vector space with quadratic form $\Q_{\ol{V}}=-\Q_V$. 
Then 
$$ \kappa:\,V\oplus \ol{V}\to W,\ \ (v,w)\mapsto
\big(v+w,\hh (v-w)\big) $$
is an isometry, with inverse $\kappa^{-1}(x,y)=(x/2+y,x/2-y)$.   
Using the isomorphism 
$$ \Cl(V)\otimes \Cl(\ol{V})=\Cl(W)$$
to view $\Cl(V)$ as a subalgebra of $\Cl(W)$, the homomorphism 
$\varrho:\,\Cl(V)\to \mf{gl}(\wedge V)$ is simply the restriction of 
$\pi$.
The inclusion $\Cl(V)\to \Cl(W)$ restricts to an inclusion of 
Spin groups, $\Spin(V)\to \Spin(W)$. The corresponding 
inclusion $ \iota:\,\SO(V)\to \SO(W)$ is given by 
$$ \iota:\,\SO(V)\to \SO(W),\ \ C\mapsto \kappa\circ \left(
\begin{array}{cc}C&0\\0&I
\end{array}\right)\circ \kappa^{-1}
=\left(\begin{array}{cc}\f{1}{2}(C+I)&{C-I}\\ \f{1}{4}(C-I)&\f{1}{2}(C+I)
\end{array}\right).
$$
\begin{proposition}\label{prop:f1}
Let $C\in\SO(V)$ with $\det(C-I)\not=0$, and suppose that 
$D\in \o(V)$ is invertible and commutes with $C$. Then 
there is a unique factorization 
\begin{equation}\label{eq:factorization}
 \iota(C)=
\left(\begin{array}{cc}
I&0\\ E_1&I
\end{array}\right)
\left(\begin{array}{cc}
I&D\\0&I
\end{array}\right)
\left(\begin{array}{cc}
I&0\\ E_2&I
\end{array}\right)
\left(\begin{array}{cc}
R&0\\0&(R^{-1})^t
\end{array}\right)
\end{equation}
such that $E_1,E_2\in\o(V)$ and $R\in \GL(V)$ commute with 
$C$ and $D$. One finds 
$$
E_1=\f{1}{2}\f{C+I}{C-I}-\f{1}{D},\ \
E_2=\f{1}{D^2}\left(\f{C-C^{-1}}{2}-D\right),\ \  R=\f{D}{I-C^{-1}}.$$
\end{proposition}
\begin{proof}
Working out the matrix product on the right hand side 
Equation \eqref{eq:factorization} reads
$$\left(\begin{array}{cc}\f{1}{2}(C+I)&{C-I}\\ \f{1}{4}(C-I)&\f{1}{2}(C+I)
\end{array}\right)
=\left(\begin{array}{cc}
(I+DE_2)R& D(R^{-1})^t\\(E_1+E_1DE_2+E_2)R&(I+DE_1)(R^{-1})^t
\end{array}\right) 
$$
It is straightforward to check this equality with the given 
formulas for $R,E_1,E_2$. Conversely, equality of the upper
right corners gives $R$, and then the diagonal entries  
give our formulas for $E_1,E_2$. 
\end{proof}

The factorization \eqref{eq:factorization} gives rise to 
a factorization for any $\hat{C}\in \Spin(V)$ 
covering $C$. The first three factors lift as in (i),(ii) 
above, and a lift $\wh{R}\in \on{ML}(V)$ of $R=D/(I-C^{-1})$ is determined 
by the choice of lift $\hat{C}$ of $C$. Using the known 
action of each factor in the representation $\pi$ we obtain:
\begin{proposition}\label{prop:main}
Suppose $\hat{C}\in\Spin(V)$ maps to $C\in \SO(V)$ with
$\det(C-I)\not=0$, and that $D\in\o(V)$ is invertible and commutes
with $C$. Let $E_1,E_2\in\o(V)$ and $\hat{R}\in \on{ML}^+(V)$ be as
above. Then the operator $\varrho(\hat{C})$ on $\wedge V$ has the
following factorization:
\begin{equation}\label{eq:rho}
\varrho(\hat{C}).\alpha=
\f{\exp(\iota_{\lambda(E_1)})\exp(\lambda(D))\exp(\iota_{\lambda(E_2)})
R.\alpha}{ {|\det|}^{1/2}(\hat{R})}.
\end{equation}
In particular, the symbol of $\hat{C}$ is given by the formula, 
\begin{equation}\label{eq:symb}
 q^{-1}(\hat{C})=
\f{\exp(\iota_{\lambda(E_1)})
\exp(\lambda(D))}{
 {|\det|}^{1/2}(\hat{R})
}.
\end{equation}
\end{proposition}
\vskip.2in

Given $\hat{C}$, there may be many natural choices of $D$ with the required 
properties, leading to different formulas for the symbol of $\hat{C}$.

\begin{proposition}
\begin{enumerate}
\item
Suppose $\hat{C}\in \Spin(V)$ maps to $C\in \SO(V)$, with  
$\det(C+I)\not=0$. Then 
\begin{equation}\label{eq:I}
{\textstyle 
q^{-1}(\hat{C})=\pm\  {\det}^{1/2}(\f{C+I}{2})\exp(2\lambda(\f{C-I}{C+I})).}
\end{equation}
\item Suppose $\hat{C}\in \on{Pin}(V)$ maps to $C\in \O(V)$, with  
$\det(C-I)\not=0$. Then
\begin{equation}\label{eq:II}
{\textstyle 
q^{-1}(\hat{C})=\pm\  {\det}^{1/2}(I-C^{-1})\exp(
\iota_{\lambda(\f{1}{2} \f{C+I}{C-I} )}) \ \d\Vol,}
\end{equation}
where $\d\Vol$ is the volume form on $V$ given by the quadratic form and some
choice of orientation. 
\end{enumerate}
\end{proposition}

\begin{proof}
Let us first assume that $\hat{C}\in \Spin(V)$ and that { both} 
$\det(C+I)\not=0$ and 
$\det(C-I)\not=0$. (In particular, $\dim V$ must be even.) 
The first formula is obtained from the choice 
\begin{equation}\label{eq:D}D=2(C-I)/(C+I),\end{equation}
since $E_1=0$ in this case. Let $D_t:=tD$ for $t>0$. As $t\to \infty$, 
$$ \exp(\lambda(D_t))=t^{\dim V/2} {\det}^{1/2}(D)\d\Vol+O(t^{\dim V/2-1}),$$
where ${\det}^{1/2}(D)$ the Pfaffian corresponding to the choice of 
orientation, and 
$$ {\textstyle 
{\det}^{1/2}(R_t)=t^{\dim V/2}{\det}^{1/2}(\f{D}{I-C^{-1}}).}
$$
The factors $t^{\dim V/2}{\det}^{1/2}(D)$ cancel, and taking the limit
$t\to \infty$ we obtain the second formula. By continuity, 
one can drop the assumption $\det(C-I)\not=0$ in the first formula 
and the assumption $\det(C+I)\not=0$ in the second formula. 
The first formula also holds if $\dim V$ is odd, by restricting  
the formula for $V\oplus\R$. Similarly, to extend the second formula 
to $\wh{C}\in \Pin(V)$ one replaces $V$ with $V\oplus\R$ and $C$ with
$\tiny(\!\!\begin{array}{cc}C\!\!\!\!\!\!&0\\0\!\!\!\!\!\!&-1
\end{array}\!\!)$. 
\end{proof}

The sign ambiguity is resolved if $C=\exp(A)$ and 
$\hat{C}=\exp(\gamma(A))$. In this case 
Equation \eqref{eq:I} reduces to 
\eqref{eq:c1}, and Equation \eqref{eq:II} to \eqref{eq:c2}.
Similarly, 
Equation \eqref{eq:c3} is obtained from \eqref{eq:II} by 
the choice $C=c\exp(A)$. 
The function $\S$ appears if we restate \eqref{eq:rho} for
$\hat{C}=\exp(\gamma(A))$ and the choice $D=A$. For the time 
being, we treat $\S$ as a meromorphic function of $A$.
\begin{proposition}\label{prop:sym1}
Suppose $A\in\o(V)$ has no eigenvalues in the set $\tpi \Z\backslash\{0\}$. 
Then the operator $\varrho(\exp(\gamma(A)))$ on $\wedge V$ has the following
factorization:
\begin{equation}\label{eq:fac2}
\varrho(\exp(\gamma(A)))
=\iota_{\S(A)}\circ \exp(\lambda(A))
\circ \exp(\iota_{\lambda(g(A))})\circ j^L(A)^{-1}.
\end{equation}
In particular, 
\begin{equation}\label{eq:fac3}
q^{-1}(\exp(\gamma(A)))=\iota_{\S(A)}\,\exp(\lambda(A))
\end{equation}
\end{proposition}
\begin{proof}
The assumption on $A$ implies that $f(A),g(A),j^L(A),j^R(A)$ 
are all well-defined and that $j^L(A)$ is invertible. 
If $\dim V$ is even, the two sides of \eqref{eq:fac2} are equal since they 
agree on the open dense subset where $A$ is invertible, 
by Proposition \ref{prop:main}. The odd-dimensional 
case follows by restricting the identity for $V\oplus \R$. 
\end{proof}
Equation \eqref{eq:fac3} gives Theorem \ref{th:factor}
for $\phi=\{0\}$. We now show that in fact, it implies the general 
case.

\begin{proof}[Proof of Theorem \ref{th:factor}]
Let $\ti{V}=V\oplus E$. Fix a non-degenerate 
quadratic form 
$\Q_E$ on $E$, and let $\ti{V}$ be equipped with the quadratic form  
$\Q_{\ti{V}}=\Q_V\oplus \eps \Q_E$ for $\eps>0$. Then
$\lambda(A)-{\sum}_a e_a \phi^a=\lambda(\ti{A})$
with
$$ \ti{A}=\left(\begin{array}{cc} A&-\eps \phi^t\\\eps\phi&0\end{array}\right)
\in\o(\ti{V}).$$
By the above, $q\circ \iota_{\S(\ti{A})}\exp(\lambda(\ti{A}))
=\exp(\gamma(\ti{A}))$. Equation \eqref{eq:key} follows 
by letting $\eps\to 0$. 
It remains to show that $\S$ is analytic everywhere. Equation 
\eqref{eq:key} says that 
$$ \iota_{\S(A)}\alpha(A)=\beta(A)$$
where $\alpha(A)=e^{\lambda(A)-\sum_a e_a\phi^a}$ and $\beta(A)=\
q^{-1}\Big(e^{\gamma(A)-\sum_a e_a{\phi}^a} \Big)$ are 
differential forms depending analytically on $A$. Take $E=V$ and
$\phi=\id$. Then $-\sum_a e_a \phi^a\in\wedge^2(V\oplus E)$ is
non-degenerate.  Equivalently, the form $\Gamma$ given as the top
degree part of $\exp(-\sum_a e_a \phi^a)$ is a volume form. Let $*$
denote the star operator defined by $\Gamma$,
i.e. $\zeta=*\iota_\zeta\Gamma$ for $\zeta\in \wedge(V\oplus
E)$. Then $\iota_{\S(A)}\alpha(A)=\beta(A)$ is equivalent to
$\S(A)\wedge *\alpha(A)=*\beta(A)$.  Since the top form degree part of
$\alpha(A)$ coincides with $\Gamma$, the constant term of $*\alpha(A)$
is equal to $1$. Hence $*\alpha(A)^{-1}$ is well-defined and depends 
analytically on $A$. We obtain
$$ \S(A)=*\beta(A)\wedge (*\alpha(A))^{-1},$$
showing explicitly that $\S(A)$ is analytic everywhere. 
\end{proof}
As another application of the factorization formula \eqref{eq:fac2}
we prove the alternative formula  \eqref{eq:alter} for the 
Clifford exponential: 
\begin{proof}[Proof of  \eqref{eq:alter}]
Apply Equation \eqref{eq:fac2}
to $\alpha=\exp(-\sum_a e_a\psi^a )\in \wedge(V)\otimes 
\wedge(E)$. The left hand side is 
$$ q^{-1} \big(\exp(\gamma(A))\exp(-\sum_a e_a\psi^a) \big).$$
To compute the right hand side, we first note that
$$ j^L(A)^{-1} \sum_a e_a\psi^a=\sum_a (j^L(A)^{-1}.e_a)\phi(j^R(A)e^a)
=\sum_a e_a \phi^a.$$
Hence $ j^L(A)^{-1}.\exp({-\sum_a e_a\psi^a})
=\exp({-\sum_a e_a \phi^a})$. Furthermore,
$$ \exp({\iota_{\lambda(g(A))}}).\exp({-\sum_a  e_a{\phi}^a})=
\exp({\varpi(A)})\ \exp({-\sum_a  e_a{\phi}^a}).$$
Hence the right hand side of \eqref{eq:fac2} becomes 
$$
\exp(\varpi(A))\ 
\iota_{\S(A)}\big(\exp(\lambda(A))\exp(-\sum_a e_a \phi^a)\big).
$$
\end{proof}

\bibliographystyle{amsplain}

\def\cprime{$'$} \def\cprime{$'$} \def\cprime{$'$}
\providecommand{\bysame}{\leavevmode\hbox to3em{\hrulefill}\thinspace}
\providecommand{\MR}{\relax\ifhmode\unskip\space\fi MR }
\providecommand{\MRhref}[2]{%
  \href{http://www.ams.org/mathscinet-getitem?mr=#1}{#2}
}
\providecommand{\href}[2]{#2}

\end{document}